\documentclass[11pt, abstract=yes]{scrartcl}
\usepackage[utf8]{inputenc}
\usepackage[T1]{fontenc}
\usepackage{amsthm,amssymb,amsmath}
\usepackage{graphicx}
\usepackage{enumerate}
\usepackage{hyperref}
\usepackage{tikz-cd}
\theoremstyle{definition}
\newtheorem{definition}{Definition}
\newtheorem{mexample}[definition]{Example}
\theoremstyle{remark}
\newtheorem{remark}[definition]{Remark}
\theoremstyle{plain}
\newtheorem{lemma}[definition]{Lemma}

\newtheorem{proposition}[definition]{Proposition}
\newtheorem{theorem}[definition]{Theorem}
\newtheorem{corollary}[definition]{Corollary}

\newcommand{\set}[1]{\left\{{#1}\right\}}
\newcommand{\vek}[1]{\boldsymbol{#1}}
\newcommand\setsuchas[2]{\left\{\,{#1}\,\vrule\,{#2}\,\right\}}
\newcommand{\NN}{{\mathbb{N}}}
\newcommand{\ZZ}{{\mathbb{Z}}}
\newcommand{\CC}{{\mathbb{C}}}
\newcommand{\RR}{{\mathbb{R}}}
\newcommand{\QQ}{{\mathbb{Q}}}
\newcommand{\Newton}{\mathrm{New}}
\newcommand{\conv}{\mathrm{conv}}

\newcommand{\restrict}{\rho}
\newcommand{\logiso}{\mathrm{log}}
\newcommand{\expiso}{\mathrm{exp}}
\newcommand{\vol}{\mathrm{vol}}
\newcommand{\mingen}{\mathrm{mg}}
\newcommand{\MA}{\mathcal{A}}
\newcommand{\amult}{\mathrm{w}}
\newcommand{\SUFM}{\mathcal{M}}
\author{Jan Snellman\thanks{Department of Mathematics, Linköping University, 58183 Linköping, Sweden}\thanks{jan.snellman@liu.se}}
\date{2025-02-23}
\title{Infinite Minkowski sums of lattice polyhedra}
\hypersetup{
 pdfauthor={Jan Snellman\thanks{Department of Mathematics, Linköping University, 58183 Linköping, Sweden}},
 pdftitle={Infinite Minkowski sums of lattice polyhedra},
 pdfkeywords={},
 pdfsubject={},
 pdfcreator={Emacs 30.1.50 (Org mode 9.7.22)}, 
 pdflang={English}}
\usepackage[style=numeric]{biblatex}
\begin{document}

\maketitle

\begin{abstract}
Artinian integrally closed monomial ideals are characterized by their Newton polyhedra, which are
lattice polyhedra inside the positive orthant having the positive orthant as their recession cone.
Multiplication of such ideals correspond to Minkowski addition of their Newton polyhedra.
In two dimensions, the isomorphic monoids of artinian, integrally closed monomial ideals under multiplication,
or the class of lattice polyhedra described above, under Minkowski addition, are free abelian, as proved by
Crispin-Quinonez.

Bayer and Stillman considered so-called \emph{monomial submodules} of the Laurent polynomial ring.
Inspired by this, we consider a family of such monomial submodules that can be (uniquely) expressed as an infinite
product of monomial submodules isomorphic to integrally closed monomial ideals. Geometrically,
their Newton polyhedras are expressed as an infinite Minkowski sum of \textbf{simple} lattice polyhedra.
This gives another example of a \textbf{topological ufm} i.e. a topological abelian monoid in which
every element can be uniquely written as a convergent (possibly infinite) product of irreducibles.
\end{abstract}
\section{Introduction}
\label{sec:orgf56db31}

Let \(I \subset \CC[\vek{x}]\) be a monomial ideal, that is, \(I\) is the \(\CC\)-span of
its monomials. Put \(J=\log(I)\),  the set of exponent vectors of
monomials in \(I\). Then \(J \subset J + \NN^{n}\), so \(J\) is a \textbf{monoid ideal}, and its
\textbf{integral closure} is
\begin{equation}
\label{eq:monid-ic}
 \bar{J} = \{\vek{\alpha} \in \ZZ^{n} \, \mid \, \exists r > 0: r \vek{\alpha} \in r J\}.
\end{equation}
It is well known \autocite{Ebud:View,V:MonAlg} that
\begin{itemize}
\item the integral closure of \(J\) is the set \(\conv(J) \cap \NN^{n}\),
the set of lattice points in the convex hull (in \(\RR^{n}\)),
\item the integral closure of \(\bar{I}\) is the \(\CC\)-span of \(\bar{J}\), hence is the
monomial ideal generated by the monomials whose exponent vectors lie in the convex hull
of \(\log(I\).
\end{itemize}

Zariski \autocite{ZariskiSamuel2,ZariskiIC} showed that the set of
\(\mathfrak{m}\)-primary, integrally closed ideals in a local, two-dimensional
domain constitutes a free abelian monoid under multiplication.
Crispin Quinonez \autocites{CrispinPhD}[][]{Crispin:ICop}[][]{CR:IM}
and Gately \autocite{Gately:OneFibered} studied integrally closed
monomial ideals in two and three variables. The main result in \autocite{Crispin:ICop} is
an explicit description of the unique factorization of \(\mathbf{m}\)-primary,
integrally closed monomial ideals in two variables into \textbf{simple} ideals.

Since \(\mathfrak{m}\)-primary integrally closed monomial ideals correspond to
can be identified with certain \textbf{lattice polyhedra} in \(\RR_+^n\) whose
complement has finite volume. Since multiplication of ideals correspond to
Minkowski addition of their polyhedra, the above result for two-dimensional
monomial ideals can be stated as follows: the monoid (under Minkowski
addition) of lattice polyhedra in \(\RR_+^2\) which
\begin{itemize}
\item are stable under Minkowski addition with \(\RR_+^2\),
\item have complement to \(\RR_+^2\) which has finite volume,
\end{itemize}
is free abelian.

In the study of resolutions of monomial ideals, so-called \textbf{monomial modules}
occur naturally \autocite{BayerSturmfels:Cellular}. A monomial module \(M\) is a
\(\CC[\vek{x}]\)-submodule of the ring \(\CC[\vek{x},\vek{x}^{-1}]\) of Laurent
polynomials which is generated by Laurent monomials \(\vek{x}^{\vek{\alpha}}\),
\(\vek{\alpha} \in \ZZ^n\). Thus \(M\) is the \(\CC\)-span of its Laurent monomials.
Put \(J = \log(M) \subset \ZZ^{n}\). Then \(J\) satisfies \(\NN^{n} + J \subset J\).
We call such a subset a \textbf{monoid (sub)module}.
Similar to the case of monoid ideals, we define the integral closure of
\(J\) by (\ref{eq:monid-ic}). It is again easy to see that the integral closure of
\(J\) is the set of lattice points in the convex hull.

By analogy with the ideal case, we call the monomial module \(M\) integrally closed
if and only if the monoid module \(J\) is; so by definition \(M\) is integrally closed
if and only if the convex hull of
\(\log(M)\) contains no additional lattice points.

Note that a finitely generated,
integrally closed monomial modules is isomorphic to a unique integrally closed
ideal. However, there are plenty of non-finitely generated monomial modules,
even if we impose the restriction that these modules should be what Bayer and
Sturmfels \cite{BayerSturmfels:Cellular} call \textbf{co-artinian}: i.e., containing no
infinite descending (w.r.t. divisibility) sequence of Laurent monomials. In
what way, may one ask, do integrally closed, co-artinian monomial modules
factor? Since the corresponding question for monomial ideals is solved for
\(n=2\) but open for larger \(n\), it is reasonable to restrict to the case
\(n=2\).

\begin{figure}[htbp]
\centering
\includegraphics[width=0.4 \textwidth]{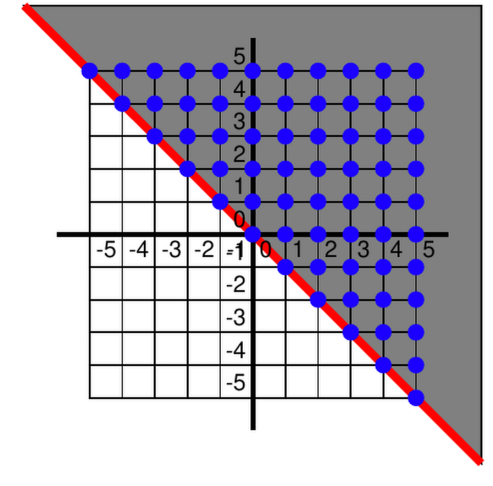}
\caption{\label{fig:infmax}An infinite power of the maximal ideal}
\end{figure}

A monomial module \(M\) is in some sense the limit of the monomial ideals
obtained by intersecting its Newton polyhedron with translates of the positive
quadrant, \(I_a(M)=x^ay^a M \cap \CC[x,y]\). If \(M\) is the co-artinian monomial
module with minimal generators \(x^ay^{-a}\), \(a \in \ZZ\), then the sequence of
monomial ideal obtained in this way consists of powers of the maximal ideal
\((x,y)\). Thus, the factorization of \(M\) into simple integrally closed
modules would, if it existed, be \(M = (x,y)^\infty\). Se Figure \ref{fig:infmax}.

The problem here is that the Newton polyhedron \(\Newton(M)\) has one
unbounded face not parallel with the \(x\) or \(y\) axis. If we restrict to
(generalized) polyhedra which have \(\RR_+^2\) as their recession cone, this particular
problem does not arise. As we shall see, if we restrict ourseleves even further,
to a very particular type of integrally closed monomial submodules of \(\CC[x,x^{-1},y,y^{-1}]\),
they can in fact be uniquely  factored into a
\emph{convergent, countable product} of (translations of ) simple integrally closed monomial ideals.
This yields another example of a \textbf{topological ufm}, which is the author's term for
abelian topological monoid in which every element can be uniquely written as a convergent, possibly infinite,
product of irreducibles \autocite{Snellman:TopUFD,Snellman:Invlimerrata,Snellman:topfact}.

\begin{mexample}
For a concrete example, consider the monomial module \(N\) with minimal
generators
\begin{displaymath}
  \set{1} \cup
  \setsuchas{x^{1+2+\cdots +r} y^{-r}}{r \in \NN_{+}},
\end{displaymath}
and let \(M\) be its integral closure. Then
\begin{displaymath}
  N \simeq \prod_{i=1}^\infty E_i,
\end{displaymath}
where \(E_i\) is the integral closure of the monomial ideal \((x^i,y)\). For
the corresponding Newton polyhedra, this becomes
\begin{displaymath}
  \Newton(N) = \sum_{i=1}^\infty \Newton(E_i),
\end{displaymath}
and infinite Minkowski sum of polyhedra.
This is illustrated in Figure \ref{fig:infmin}.
\begin{figure}[htbp]
\centering
\includegraphics[width=0.6 \textwidth]{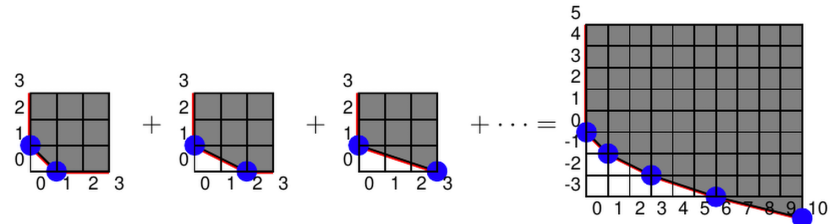}
\caption{\label{fig:infmin}The Newton polyhedra of a non-finitely generated integrally closed monomial module, expressed as an infinite Minkowski sum of simple polyhedra.}
\end{figure}
\label{ex:concrete}
\end{mexample}

\begin{remark}
While co-artinian monomial modules \(M \subset \CC[\vek{x}]\) are nominally the object of study in this article,
we will in actual fact only deal with the following corresponding objects:
\begin{itemize}
\item the monoid module \(J \subset \ZZ^{n}\),
\item the Newton polyhedron \(\Newton(M) \subset \RR^{n}\), and
\item the formal power series in \(\ZZ[[\vek{x}]]\) which is the generating function of \(J\) (or the \emph{integer transform} of
the Newton polyhedron)
\end{itemize}
Furthermore, we will study only the case of \(n=2\), and restrict ourselves to some very special \(J\).
The purpose of the article is to present the simplest possible example of unique infinite factorization
(of lattice polyhedra or of monoid modules). So while the inspiration comes from commutative algebra,
the results and methods are from combinatorics, semigroup theory, and convex geometry.
\end{remark}
\section{Monomial modules and monoid submodules}
\label{sec:org31aeb97}
The additive monoid \(\NN^{2}\) is a submonoid of \(\ZZ^{2}\), its \emph{difference group}
(sometimes called \emph{Grothendieck group}). This gives an inclusion of semigroup
algebras \(\CC[\NN^{2}] \subset \CC[\ZZ^{2}]\), which gives the \emph{Laurent polynomial ring}
\(\CC[\ZZ^{2}]\) a natural \(\CC[\NN^{2}]\) module structure.

Note that \(\CC[\NN^{2}] = \CC[x,y]\) and that \(\CC[\ZZ^{2}]=\CC[x,x^{-1},y,y^{-1}]\).
We will write \(T=\setsuchas{x^{a}y^{b}}{(a,b) \in \ZZ^{2}}\) and \(S=\setsuchas{x^{a}y^{b}}{(a,b) \in \NN^{2}}\)
for the corresponding, isomorphic but multiplicatively written, monoids, and use
\(\logiso: T \to \ZZ^{2}\) and \(\expiso: \ZZ^{2} \to T\) for the natural isomorphisms.
We will also view \(T\) as a subset (or indeed a \(\CC\)-basis) of \(\CC[\ZZ^{2}]\) and
\(S\) as a subset (or indeed a \(\CC\)-basis) of \(\CC[\NN^{2}]\), and alternatively
write these semigroup rings as \(\CC[T]=\CC[\ZZ^{2}]\) and \(\CC[S]=\CC[\NN^{2}]\).

A \emph{monomial ideal} \(I \subset \CC[S^{}]\) is an ideal generated by its monomials \(I \cap S\);
thus \(I\) is the \(\CC\)-vector space span of \(I \cap S\).
The corresponding \emph{monoid ideal} \((I \cap S) \subset S\) is a filter with respect to
the divisibility order on \(S\); similarly \(\logiso(I \cap S) \subset \NN^{2}\) is a filter
with respect to the cartesian order \(\le\) on \(\NN^{2}\).
An equivalent definition of a monoid ideal \(J \subset S\) is to demand that \(SJ \subset J\).

The \emph{integral closure} \(\bar{I}\) of \(I\) is the set of all elements  \(u \in \CC[S]\)
such that
\begin{equation}
\label{eq-ic}
u^{m} + a_{1} u^{m-1} + \cdots + a_{m-1}u + u_{m} = 0, \qquad a_{i} \in I^{i}.
\end{equation}
This is a monomial ideal \autocite{V:MonAlg} so the \(\CC\)-span of \(\bar{I} \cap S\); the monoid ideal
\(\bar{J}=\bar{I} \cap S\) is the integral closure of \(J=I \cap S\), where the integral closure of a monoid ideal
is
\begin{equation}
\label{eq:icmonid}
    \log(\bar{J}) = \setsuchas{\vek{\alpha} \in \NN^2}{\exists r>0: \, r  \vek{\alpha} \in r \log(J)}.
\end{equation}

It holds that
\begin{equation}
\label{eq:jbar}
    \logiso(\bar{J}) = \conv(\logiso(J)) \cap \NN^{2}
\end{equation}
where the convex hull is taken inside \(\RR^{2}\). The ideal \(I\), (or the monoid ideal \(J\)), is \emph{closed}
whenever it is equal to its closure.

The so-called \emph{Newton polyhedra}  of a monomial ideal
\(I \subset \CC[S]\) (or of its associated monoid ideal \(J\), or of \(\log(J)\)) is
\[\Newton(I) = \conv(\log(I \cap S)) \subset \RR^{2};\]
it is in fact contained in the positive quadrant \(\RR^{2}_{{+}}\).
Note that \(I\) and \(\bar{I}\) have the same Newton polyhedra.
If \(I\) is artinian
(equivalently, \(\mathfrak{m}\)-primary) then \(S\setminus (I \cap S)\) is finite, and
\(\vol(\RR_{+}^{2} \setminus \Newton(I)) < \infty\),
and the recession cone of \(\Newton(I)\) is \(\RR_{+}^{2}\). For the latter condition to hold,
it is enough that \(I\) contains some monomial which is not a power of \(x\) and some monomial which
is not a power of \(y\).

For \(P,Q\) two polyhedra in \(\RR^{2}\), their
\emph{Minkowski sum} is
\[P + Q = \setsuchas{\vek{a} + \vek{b}}{\vek{a} \in P, \vek{b} \in Q}.\]
It holds that
\[\Newton(I_{1}I_{2}) = \Newton(I_{1}) + \Newton(I_{2})\]
whenever \(I_{1},I_{2}\) are integrally closed monomial ideals.

When \(J \subset S\) is a monoid ideal, the set \(\mingen(J)\) of its minimal elements with respect to divisibility
is finite, and
\[J = \cup_{m \in \mingen(J)}mS.\]
Similarly, we define the minimal generators of a monomial ideal \(I\) as the monomials corresponding
to the minimal generators of \(J = I \cap S\); they do indeed generate \(I\) minimally as an ideal of \(\CC[S]\).
When \(I\) is integrally closed, its minimal generators corresponds to the vertices of the Newton polyhedron.

Bayer and Stillman \autocite{BayerSturmfels:Cellular} called a \(\CC[S]\)-submodule of \(\CC[T]\) a \emph{monomial module}
and used such objects in their study of homological properties of monomial ideals.
If \(M \subset \CC[T]\) is such a monomial module, then it is the \(\CC\)-span of \(L = M \cap T\).
The subset \(L \subset T\) is not a monoid ideal of \(T\), nor a submonoid; it is however a filter with respect to
divisibility, and satisfies \(SL \subset L\) (but \textbf{not} \(TL \subset L\)).
We say that \(L\) is an \(S\)-submodule of the free abelian group \(T\).
We define the integral closure of \(L\) as
\begin{equation}
\label{eq:Lbar}
    \bar{L} = \setsuchas{\vek{\alpha} \in \ZZ^2}{\exists r>0: \, r  \vek{\alpha} \in r J}
\end{equation}
It holds, similar to the monoid ideal case, that
\begin{equation}
\label{eq:Lbarconv}
    \logiso(\bar{L}) = \conv(\logiso(L)) \cap \ZZ^{2}.
\end{equation}
We will in this article  refer to
\(\bar{M}\), the \(\CC\)-span of \(\bar{L}\), as the integral closure of \(M\), and call
\[\conv(\logiso(L)) \subset \RR^{2}\]
the Newton polyhedron of \(M\), even though it strictly speaking need not
be a polyhedron, since it can be the intersection of infinitely many halfplanes.

\begin{figure}[htbp]
\centering
\includegraphics[width=0.4 \textwidth]{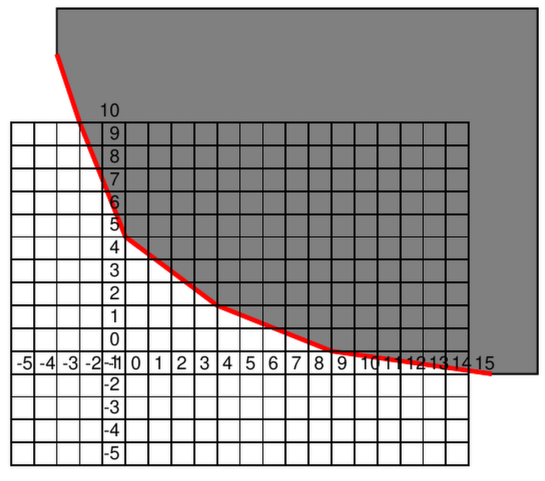}
\caption{\label{fig:coartinian}A co-artinian monomial module}
\end{figure}

Following Bayer and Stillman, we call a monomial module \(M\)  \emph{co-artinian}
if
\[L = M \cap T \subset T\]
contains no infinite descending (with respect to divisibility) chain;
equivalently, if \(L\) contains no infinite principal order ideal. A co-artinian monomial module can still
extend infinitely far ``to the left'', as Figure \ref{fig:coartinian} shows.

\begin{remark}
By a slight abuse of notation, we call \(L\subset T\), or \(\log(M) \subset \ZZ^{2}\) co-artinian
whenever \(M\) is co-artinian.
\end{remark}
\section{Factorization of integrally closed monomial ideals}
\label{sec:orge569d4f}
The following are known for monomial ideals in two variables and three variables: \autocite{CR:IM,V:MonAlg,Crispin:ICop,CrispinPhD,Gately:OneFibered}:
\begin{itemize}
\item the integral closure of a monomial ideal is a monomial ideal,
\item the product of integrally closed monomial ideals is an integrally closed monomial ideal,
\item multiplication of integrally closed monomial ideals correspond to \emph{Minkowski addition} of their
\emph{Newton polyhedra},
\item in two variables, the monoid of integrally closed artinian monomial ideals under ideal multiplication is
a free abelian monoid of infinite rank, with basis given by the \emph{simple ideals}
\(E_{r,s} =(x^{r}, y^{s})\), \(\gcd(r,s)=1\),
\item in three variables, factorization is no longer unique, nor is the length of
factorizations of a given ideal always the same.

In Figure \ref{fig:icfact} we illustrate the factorization
\[
  E_{4,3} E_{5,2} = (x^{4},y^{3})(x^{2},y^{2}) = (x^{9}, x^{4}y^{2}, y^{5}).
  \]
\begin{figure}[htbp]
\centering
\includegraphics[width=0.6 \textwidth]{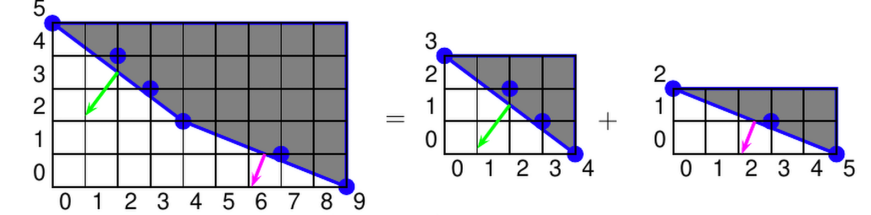}
\caption{\label{fig:icfact}Factorization of \(I=(x^{9}, x^{4}y^{2}, y^{5})\)}
\end{figure}

The corresponding monoid ideals, Newton polyhedron, and the boundary of the Newton polyhedron,
is shown in Figure \ref{fig:idnebo}.
\begin{figure}[htbp]
\centering
\includegraphics[width=0.6 \textwidth]{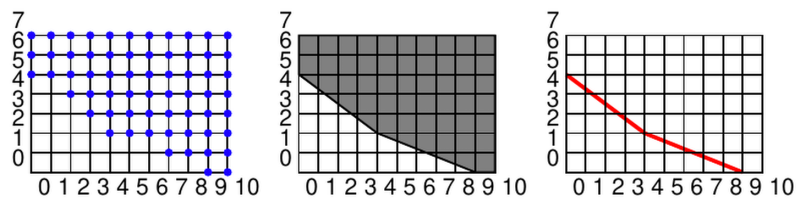}
\caption{\label{fig:idnebo}Associated monoid ideal and Newton polyhedra of \(I=(x^{9}, x^{4}y^{2}, y^{5})\)}
\end{figure}
\end{itemize}

\begin{remark}
In what follows, we will freely switch between the additive monoid submodule \(L \subset \ZZ^{2}\)
and the isomorphic multiplicative monoid submodule \(\exp(L) \subset T\),
sometimes calling a sum a ``factorization''. We trust that the reader will not be confused.
\end{remark}

We say that two monoid submodules \(L_{1}, L_{2} \subset \ZZ^{2}\) are equivalent if \(L_{2}= \vek{\alpha} + L_{1}\) for
some \(\vek{\alpha} \in \ZZ^{2}\). If \(L_{1},L_{2}\) are integrally closed, this is happens iff there
is an integral vector translating one Newton polyhedron into the other. A finitely generated
monoid module \(L\) is thus equivalent with a monoid ideal of \(\NN^{2}\).
Furthermore, if \(L\) is integrally closed and \(\Newton(L)\)
has recession cone \(\RR_{+}^{2}\), we can translate it so that it ``aligns'' with
\(\RR_{+}^{2}\). In this case, \(L\) is equivalent to the complement of a finite order ideal
in \(\NN^{2}\) (or to an artinian integrally closed monomial ideal).
\section{Infinite factorization}
\label{sec:org51df351}
We will be using some facts and notions about topological spaces, topological monoids, and topological groups, for which references may be found in \autocite{Bourbaki:Topology,Bourbaki:Ensemble,Kelley:Topology}.
The author studied \autocite{Snellman:TopUFD,Snellman:Invlimerrata,Snellman:topfact} inverse limits
of (discretely topologized) unique factorization domains.
The resulting topological rings, called \textbf{topological unique factorization domains},
have the property that any element can be uniquely factored
as a \textbf{convergent} (possibly infinite) product of irreducible elements. Examples of such
topological ufd's are
Halter-Kochs \autocites{HK:FinGenMon}[][]{HK:genpolalg} ``formal polynomials'' are \(\varprojlim k[x_{1},\dots,x_{n}]\)
  where the truncations are the ring homomorphisms setting the last variable to zero.

For a topological domain \(R\), the set of invertible elements \(R \setminus \set{0}\)
is a topological monoid, which captures the factorization porperties of \(R\).
Indeed, \(R\) is a ufd iff \(R \setminus \set{0}\) is a free abelian monoid.
Thus we call a (abelian, topological) monoid with unique convergent factorization into irreducibles,
such as the set of invertibles elements of a topological UFD, a \textbf{topological unique factorization monoid}.

The purpose of this article is to exhibit yet another interesting example of a topological unique factorization monoid.
\begin{definition}
Let
\[\SUFM = \NN^{\NN_{+}} =\setsuchas{(a_{j})_{j=1}^{\infty}}{a_{j} \in \NN}\]
 be the monoid of all sequences of non-negative
integers, with component-wise addition, and endowed with the product topology, where \(\NN\) is
discretely topologized.
\end{definition}
This means that \(\vek{a}^{j} \to \vek{a} \in \SUFM\) if and only if there exists
some function \(v: \NN_{+} \to \NN_{+}\) with the property that for each \(i\) and each \(j>v(i)\)
it holds that \(a^{j}_{i} = a_{i}\).
Similarly, the sum
\(\sum_{j=1}^{\infty}\vek{a}^{j}\)
 is convergent, with sum \(\vek{a} \in \SUFM\), if and only if for all \(i\)
 \(\vek{a}^{j}_{i}=0\) for almost all \(j\).

\begin{proposition}
Addition and inversion in \(\SUFM\) are continuous, i.e., \(\SUFM\) is a topological monoid.
Furthermore, \(\SUFM\) is a topological unique factorization monoid,
and it is the inverse limit of the surjective inverse system
\begin{equation}
\label{invsys1}
 \SUFM_{1} \leftarrow \SUFM_{2} \leftarrow \SUFM_{3} \leftarrow \cdots
\end{equation}
where each \[\SUFM_{k} = \setsuchas{(a_{j})_{j=1}^{\infty}}{a_{j}=0 \text{ for } j> k}\] is a discrete,
finite rank free abelian monoid, and the restriction maps are
\[(a_{1},\dots,a_{k},a_{k+1},0,0,\dots) \mapsto (a_{1},\dots,a_{k},0,0,\dots).\]
\end{proposition}
\begin{proof}
 The first assertion is trivial, as is the third. For the second,
 define \(\vek{e}_{j} \in \SUFM_{j} \subset \SUFM\) to be the sequence with 1 in
 the \(j\)'th coordinate, and zeroes elsewhere.
 Then all \(\vek{e}_{j}\) are irreducible, and no other elements in \(\SUFM\) are.
It follows that
 \[(a_{1},a_{2},\dots) = \sum_{j=1}^{\infty}a_{j}\vek{e}_{j}\]
 is the unique factorization of \(\vek{a}=(a_{1},a_{2},\dots)\).
\end{proof}

We now return to monoid submodules of \(\ZZ^{2}\), and define a minimal subclass such that
unique infinite factorization is guaranteed. In fact, the resulting topological monoid
will be exactly \(\SUFM\).
\begin{definition}
Define \(\MA\) to be the collection of subsets \(L \subset \ZZ^{2}\),
such that
\begin{enumerate}
\item \(L\) is a monoid submodule,
\item \(L\) is co-artinian,
\item \(L\) is integrally closed,
\item \(L\) contains \(\NN\) but no lattice point strictly smaller than \(\vek{0}\),
i.e., nothing in the third and fourth quadrant, (except \(\vek{0}\)),
\item \(P=\Newton(L)\) contains no infinite ray \(\setsuchas{t (a,b)}{t \ge 0}\) with \(b < 0 < a\),
\item There is no infinite facet of \(P\) except the positive \(y\)-axis and (possibly) an infinite
ray parallel with the \(x\)-axis.
\item The finite facets are line segments joining vertices \((c,d)\) and \((c+k,d-1)\)
with \(d < 0 < c\), i.e., they have outward normals \((-1,-k)\).
\end{enumerate}

Denote by \(\MA_{n}\) the subcollection consisting of
those \(L \in \MA\) where the outward normals of
\(\Newton(L)\) belong to \(\setsuchas{(-1,-k)}{1 \le k \le n}\).
\end{definition}
\begin{mexample}
An example of an acceptable monoid module is the right-most part of Figure \ref{fig:infmin}.
\end{mexample}

\begin{lemma}
If \(L \in \MA_{n}\) then \(P = \Newton(L)\) has an infinite face with outward normal \((0,-1)\),
and is contained in the strip \(\setsuchas{(a,b)}{b \ge -c}\) for some non-negative integer \(c\).
\end{lemma}
\begin{proof}
Obvious.
\end{proof}

\begin{theorem}
For any positive integer \(n\), the set \(\MA_{n}\) is (under addition of monoid submodules)
a free abelian module of rank \(n\), and thus isomorphic to \(\SUFM_{n}\).
\label{thm-finfac}
\end{theorem}
\begin{proof}
Take \(L \in \MA_{n}\). Then there is a unique non-negative integer \(c\) such that \((0,c) + L\)
is equivalent to the set of exponent vectors of monomials in
an artinian, integrally closed ideal \(I \subset \CC[S]\). It follows from Crispin-Quinonez characterization
of such ideals \autocite{Crispin:ICop,CR:IM,CrispinPhD} that \(I = \prod_{j} E_{(r_{j},s_{j})}\).
This factorization is unique, and the factors can be read off from the outward normals; thus
\(r_{j} \in \set{1,2,\dots,n}\) and \(s_{j}=1\). We get that
\begin{equation}
\label{eq:factranI}
I = \prod_{j=1}^{n} E_{(j,1)}^{a_{j}}, \qquad a_{j} \in \NN.
\end{equation}
In other words, \(I\) belongs to the submonoid
generated by \(E_{(1,1)}, E_{(2,1)},\dots,E_{(n,1)}\); this is a submonoid of a free abelian monoid generated
by \(n\) different irreducible elements, hence is isomorphic to a free abelian module with these generators
as a basis.

To each \(E_{(j,1)}\) with \(1 \le j \le n\) we associate the monoid module
\(F_j = (0,-1) + \log(E_{(j,1)})\).
It is clear that
\(F_j \in \MA_{j} \subset \MA_{n}\). Hence our original \(L \in \MA_{n}\) satisfies
\begin{equation}
\label{eq:factranI}
M = \sum_{j=1}^{n} F_j^{a_{j}}, \qquad a_{j} \in \NN.
\end{equation}
\end{proof}

\begin{mexample}
The monomial module \(M\) minimally generated by \(\set{1, xy^{-1}, x^{3}y^{-2}, x^{6}y^{-3}}\)
corresponds to the monoid module \(L=\set{(0,0),(1,-1),(3,-2),(6,-3)}\), which
factors as
\[L = F_{1} + F_{2} + F_{3},\]
as can be seen in Figure \ref{fig:infmin} (consider the first three terms of the infinite factorization
illustrated in the figure).
\end{mexample}

\begin{definition}
If \(L \in \MA\) we put \(\chi_{L} = \sum_{\vek{\alpha} \in L} \vek{x} \in \ZZ[ [T] ]\).
Give \(\ZZ\) the discrete topology and \(\ZZ[ [T] ]\) the product topology,
identifying \(\ZZ [ [T] ]\) with \(\ZZ^{T}\). Finally, give \(\MA\) the subspace topology via
the embedding
\[
\MA \ni L \mapsto \chi(L) \in \ZZ[ [T]].
\]
\end{definition}

\begin{proposition}
A sequence \((L_{j})_{j=1}^{\infty}\) of monoid submodules
in \(\MA\) converges to \(L \in \MA\)  if and only if
\begin{equation}
\label{eq:convprop}
\forall \vek{\alpha} \in \ZZ^{2}: \, \exists v(\vek{\alpha}) \in \NN_{+}: \, \vek{\alpha} \in L \, \iff \,
\forall \ell > v(\vek{t}): \vek{\alpha} \in L_{\ell}.
\end{equation}
Furthermore, the topology induced by \(\MA\) turns it (or \(L\)) into a topological monoid, and
each \(\MA_{n} \subset \MA\), when equipped with the subspace topology, is discrete.
\label{prop:MAtop}
\end{proposition}
\begin{proof}
Since \(\MA\) has the product topology, \(L_{j} \to L\) if and only if for each
\[[\vek{x}^{\vek{\alpha}}] \chi(N_{j}) \to [\vek{x}^{\vek{\alpha}}]  \chi(N),\]
where for \(f \in \ZZ[ [T] ]\), \([\vek{x}^{\vek{\alpha}}](f)\)
denotes the coefficient of \(\vek{x}^{\vek{\alpha}}]\) in \(f\).
However, \(\ZZ\) has the discrete topology, so
\[ [\vek{x}^{\vek{\alpha}}] \chi(N_{j}) \to [\vek{x}^{\vek{\alpha}}]\chi(N)\]
if and only if this sequence is eventually constant.
\end{proof}

It is clear that the above proposition holds when ``sequence'' is replaced by ``net''.
Thus, we get:
\begin{proposition}
Let \(A \subset \MA\), \(N \in \MA\). Then
\[
\sum_{M \in A} M = N
\]
if and only if, for every \(\vek{t} \in T\) there is some finite subset \(B \subset A\)
such that, for any finite \(B \subset C \subset A\),
\(\vek{\alpha} \in N\) if and only if \(\vek{\alpha} \in \sum_{K \in C} K\).
\label{prop:convsum}
\end{proposition}

\begin{definition}
For \(M \in \MA\), \(n\) a positive integer, let \(M(n)\) denote the
monoid module consisting of all \((a,b)\) with \(b \ge -n\).
Let \(\restrict_{n}(M) =M(c(n))\) for the smallest \(c(n)\) such that \(\Newton(M(c(n)))\) contains no
boundary segment with outward normal \((-k,-1)\) with \(k >n\).
Then \(\restrict_{n}(M) \in \MA_{n}\) and hence
\[\restrict_{n}(M) = \sum_{j=1}^{n} \amult_{M}(j,n)F_{j}\]
for some non-negative integers \(\amult_{M}(j,n)\).
\label{def:multconst}
\end{definition}

\begin{lemma}
For \(M \in \MA\) it holds that \(M(n) \to M\) as \(n \to \infty\). Furthermore, as topological spaces,
\(\MA \simeq \varprojlim_{n} \MA(n)\)
where \(\MA(n)\) is the subspace of submonoids with no lattice point \((a,b)\) with \(b < -n\).
\end{lemma}
\begin{proof}
Immediate, but note that \(\MA(n)\) is not a submonoid of \(\MA\).
\end{proof}

\begin{lemma}
Let \(c,d\) be positive integers, and let \(\vek{\alpha}=(c,-d)\). Then there are finitely many conditions
on \(a_{1},a_{2},a_{3},\dots,a_{c}\) that determines whether
\begin{equation}
\label{eq:fin-cond}
\vek{\alpha} \in \sum_{j=1}^{n} a_{j}F_{j} = L,
\end{equation}
when \(n \ge c\).
\label{lemma:finc}
\end{lemma}
\begin{proof}
First, note that if \(b_{j} \ge a_{j}\) for all \(j\) then \(\sum_{j} a_{j}F_{j} \subset \sum_{j} b_{j}F_{j}\).
We will use this to find conditions for a more restricted situation, then analyze the appropriate
relaxations.

Let \(S_{c,d} \subset \NN^{\NN_{>0}}\) be the set of integer tuples such that
\begin{enumerate}
\item For \(j>c\), \(b_{j}=0\),
\item \(\sum_{j=1}^{c}b_{j} =d\),
\item \(\sum_{j=1}^{c}jb_{j} \le c\).
\end{enumerate}
Then, using our observation, we get that
\[
(c,-d) \in \sum_{j=1}^{n}a_{j}F_{j} \quad \iff \quad \exists \vek{b} \in S_{c,d}: \vek{a} \ge \vek{b},
\]
where the ordering is the product (componentwise) order. But \(S_{c,d}\)
is in fact a subset of the set of all (integer) partitions of \(d\), since
if \((b_{1},\dots,b_{i},b_{i+1},\dots,b_{c}) \in S_{c,d}\) then \((b_{1},\dots,1+b_{i},-1+b_{i+1},\dots,b_{c}) \in S_{c,d}\). Hence
\(S_{c,d}\) is finite, and the conclusion follows.
\end{proof}

\begin{mexample}
Let \((c,d) = (6,3)\), so \(\vek{\alpha}=(6,-3)\). Then
\[
  S_{c,d} = \{(3,0,0),(2,1,0),(1,2,0),(0,3,0),(1,1,1)\},
\]
and we conclude that
\(\vek{\alpha} \in \sum_{j=1}^{n} a_{j}F_{j}\) if and only if
\begin{multline*}
  (a_{1} \ge 3) \bigvee  ((a_{1} \ge 2) \wedge  (a_{2} \ge 1)) \bigvee
  ((a_{1} \ge 1) \wedge (a_{2} \ge 2))  \\
\bigvee
  ( a_{2} \ge 3) \bigvee
  ((a_{1} \ge 1) \wedge (a_{2} \ge 1) \wedge (a_{3} \ge 1))
\end{multline*}
This is illustrated in Figure \ref{fig:arrowtrain}.

\begin{figure}[htbp]
\centering
\includegraphics[width=0.3 \textwidth]{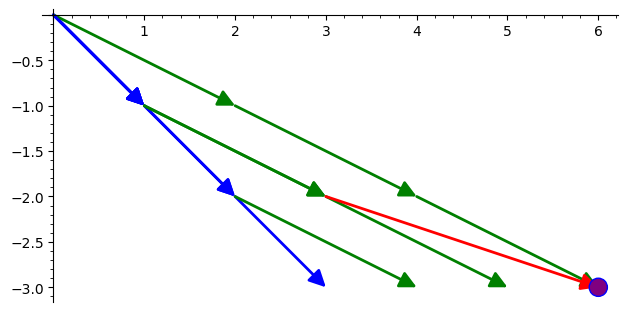}
\caption{\label{fig:arrowtrain}Ways to reach \((a,-3)\) with \(a \le 6\).}
\end{figure}
\end{mexample}

\begin{corollary}
Whenever \(j \le n,s\), it holds \(\amult_{M}(j,n) = \amult_{M}(j,s))\).
\end{corollary}
\begin{proof}
Follows from Lemma \ref{lemma:finc}.
\end{proof}

\begin{definition}
We let \(\amult_{M}(j)\) be the common value of \(\amult_{M}(j,n)\) for \(n \ge j\).
\end{definition}

\begin{theorem}
Let \((a_{j})_{j=1}^{\infty} \in \SUFM\). Then the sum
\begin{equation}
\label{eq-aseq}
\sum_{j=1}^{\infty} a_{j}F_{j} = L
\end{equation}
is convergent, with \(L \in \MA\),
and \(\amult_{L}(j) = a_{j}\) for all \(j\).

Furthermore, the map
\begin{equation}
\label{eq:topiso}
\begin{split}
\eta: \SUFM & \to \MA \\
\eta((a_{j})_{j=1}^{\infty}) & = \sum_{j=1}^{\infty} a_{j}F_{j}
\end{split}
\end{equation}
is an isomorphism of topological monoids, which for any positive \(n\)
restricts to an isomorphism of discrete, free abelian monoids of finite rank
\begin{equation}
\label{eq:disciso}
\begin{split}
\eta: \SUFM_{n} & \to \MA_{n} \\
\eta((a_{j})_{j=1}^{\infty}) & = \sum_{j=1}^{\infty} a_{j}F_{j} = \sum_{j=1}^{n} a_{j}F_{j}
\end{split}
\end{equation}
\label{thm:a-conv}
\end{theorem}
\begin{proof}
Pick a lattice point \(\vek{\alpha} = (c,-d)\) in the fourth quadrant. By Lemma \ref{lemma:finc}
there are finitely many conditions on \(a_{1},\dots,a_{c}\) which determines whether
\(\vek{\alpha} \in \sum_{j=1}^{n} a_{j}F_{j}\), irrespective of \(n\), as long as \(n \ge c\). Hence, for \(n_{1},n_{2} \ge c\),
\(\vek{\alpha}\) belongs to both, or neither, of
\(\sum_{j=1}^{n_{1}} a_{j}F_{j}\) and \(\sum_{j=1}^{n_{2}} a_{j}F_{j}\). This holds for any \(\vek{\alpha}\), and so (\ref{eq-aseq}) is convergent.
From Lemma \ref{lemma:finc} it follows that
\[
\restrict[n](L) = \sum_{j=1}^{n} a_{j}F_{j}
\]
and that \(\amult_{L}(j) = a_{j}\).
It is immediate that \(\eta(\vek{a} + \vek{b}) = \eta(\vek{a}) + \eta(\vek{b})\), and that \(\eta\) is injective.
Furthermore, every
\(L \in \MA\) can be written as the convergent sum (\ref{eq-aseq}), by ``reading off'' the normals of
\(\Newton(L\).
It is also clear that \(a_{j} = \amult_{L}(j)\) for all \(j\), so the inverse of \(\eta\) is
\[
\MA \ni L \mapsto (\amult_{L}(j))_{j=1}^{\infty.}
\]

Applying Lemma \ref{lemma:finc} we see that \(\SUFM\) and \(\MA\) are homeomorphic via \(\eta\).
\end{proof}

\begin{theorem}
We have the following homomorphisms of topological abelian monoids:
\begin{equation}
\label{eq:invsys}
\begin{tikzcd}
\SUFM \arrow[dr, two heads] \arrow[drr, two heads]
\arrow[drrr, two heads] \arrow[ddd, bend right, "\simeq"] \\
& \SUFM_1 \arrow[d, "\simeq"]& \SUFM_2 \arrow[d, "\simeq"] \arrow[l, two heads]
& \SUFM_3 \arrow[d, "\simeq"]\arrow[l, two heads] & \cdots \arrow[l, two heads] \\
& \MA_1 & \MA_2 \arrow[l, two heads] & \MA_3 \arrow[l, two heads] & \cdots \arrow[l, two heads] \\
\MA \arrow[ur, two heads] \arrow[urr, two heads] \arrow[urrr, two heads]
\end{tikzcd}
\end{equation}
Hence,
\begin{equation}
\label{eq:invlim}
\varprojlim_{n=1}^{\infty} \MA_{n} \simeq \MA \simeq \SUFM \simeq \varprojlim_{{n=1}^{\infty}} \SUFM_{n}
\end{equation}
as topological monoids.
\label{thm:is-invlim}
\end{theorem}

\begin{figure}[htbp]
\centering
\includegraphics[width=0.6 \textwidth]{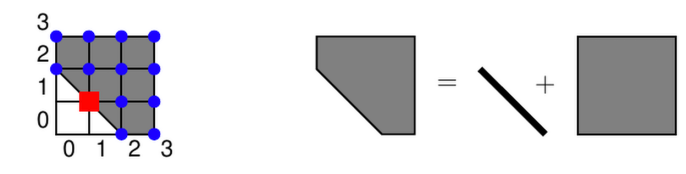}
\caption{\label{fig:MinSumLin}\((2,2) = (0,2)+N(2F_{1})\) as the Minkowski sum of line segment and \(\RR^{2}_{+}\)}
\end{figure}

\begin{corollary}
Let \(P=\Newton(M)\) for \(M \in \MA\).
Then \(P\) is in a unique way expressable as a convergent Minkowski sum of
the generalized lattice polyhedra
\[N(F_{j}) = \RR_{+}^{2} + \ell_{j,}\] where \(\ell_{j}\)
is the line segment from the origin to the lattice point \((j,-1)\).
\end{corollary}
In Figure \ref{fig:MinSumLin} one such \(N(F_{j})\) is shown.

\begin{mexample}
The first partial sums of the generalized polyhedron
\(\sum_{j=1}^{\infty} N(jF_{j})\) are shown in figure \ref{fig:infkowski}.

\begin{figure}[htbp]
\centering
\includegraphics[width=0.5 \textwidth]{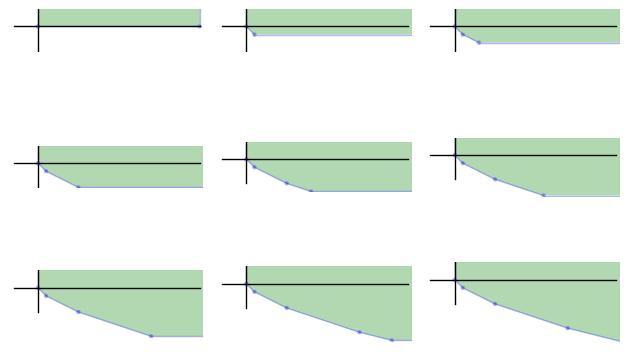}
\caption{\label{fig:infkowski}Infinite Minkowski sum of generalized lattice polyhedra.}
\end{figure}
\end{mexample}
\section{Further questions}
\label{sec:org7fe847d}
\begin{figure}[htbp]
\centering
\includegraphics[width=0.5 \textwidth]{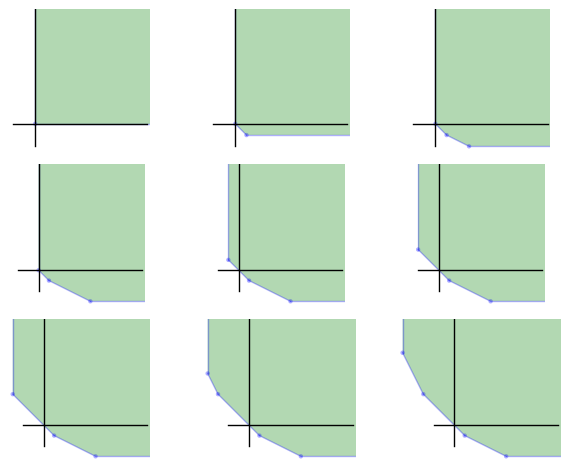}
\caption{\label{fig:anchor}Infinite Minkowski sum of generalized lattice polyhedra, properly translated}
\end{figure}

\begin{enumerate}
\item This manuscript considers a very restricted set of monomial modules.
At the very least, one should be able to
easily treat products of (translates of) \(E_{(k,1)}\) and \(E_{(1,k)}\), by anchoring the vertex where
the slope of the boundary of the Newton polyhedra increases past 1 to the origin; see figure \ref{fig:anchor}.
Some ingenuity would be needed to properly translate a general \(E_{p,q}\). Recall that
the monoid of integrally closed, \((x,y)\)-primary monomial ideals is the free abelian monoid
with these as a base. It would be satisfying to be able to study monomial modules of the form
\[M=\sum_{j=1}^{\infty} \gamma_{j} + E_{u_{j},v_{j}}\]
where \(M\) is co-artinian and
\[
   \ZZ \ni j \mapsto (u_{j},v_{j}) \mapsto u_{j}/v_{j} \in \QQ_{+}
   \]
is an order-preserving map \(\ZZ \to \QQ_{+}\). Such an \(M\) should be the product
of the corresponding simple monomid modules, \textbf{suitably translated}.
\item It is probably not possible to assign meaning to someting like \(\prod_{{q \in \QQ}} F_{q}\).
\item Infinite exponents, such as in \(M=(x,y)^{\infty}\), shown in Figure \ref{fig:infmax}, pose a problem that
may or may not be surmountable.
\item In higher dimensions factorization is no longer unique \autocite{Gately:OneFibered}.
\item Lewis \autocite{Lewis02012025} recently studied factorization of integrally closed monomial ideals
via Newton Polyhedra and the so-called \emph{integral polytope group}, also studied by Funke \autocite{FunkeIntegralPolytopeGrp}. An important
aspect of those works is factoring out with the equivalence relation
stemming from translations. We have attached the monomial modules we considered to the origin in
a particular way, in effect choosing representatives; it might be more convenient
to deal with the equivalence relation instead.
\end{enumerate}

\section{Bibliography}
\label{sec:org5ee62fa}
\printbibliography
\end{document}